\documentclass[11pt]{article} 
\usepackage[utf8]{inputenc} 
\usepackage[polutonikogreek,english]{babel}
\usepackage[or]{teubner}
\usepackage{xpatch}
\xapptocmd{\greektext}{\edef~{\string~}}{}{}
\usepackage{amsfonts}
\usepackage{amsmath}
\usepackage{amsthm}
\usepackage{titlesec}

\titleformat{\section}  
  {\fontsize{13}{15}\bfseries} 
  {\thesection.} 
  {1em} 
  {} 
  [] 

\titleformat{\subsection}  
  {\fontsize{12}{14}\bfseries} 
  {\thesubsection.} 
  {1em} 
  {} 
  [] 

\titleformat{\subsubsection}  
  {\fontsize{11}{13}\bfseries} 
  {\thesubsubsection.} 
  {1em} 
  {} 
  [] 

       \usepackage[polutonikogreek,english]{babel}

\usepackage{hyperref}
\usepackage{geometry} 
\usepackage{graphicx}
\graphicspath{ {c:/users/fpc-dprotopapas/documents/latex/} }

\usepackage{booktabs} 
\usepackage{array}
\usepackage{paralist} 
\usepackage{verbatim} 
\usepackage{subfig}
\usepackage{multirow}

\usepackage{fancyhdr} 
\usepackage[nottoc,notlof,notlot]{tocbibind} 
\usepackage[titles,subfigure]{tocloft}

\usepackage{epigraph}
\title{Simone Weil, André Weil, Bourbaki 
\\ and Pythagorean mathematics}
 \author{Athanase Papadopoulos}

\date{\today}

\begin{document}

     \maketitle
  
 \epigraph{\itshape Quoi de plus proche de nous que la Grèce ?
    Elle est plus proche de nous que nous-mêmes.
    Il est douteux que nous ayons une seule
idée importante qui n’ait été clairement conçue par les Grecs.
    (Simone Weil, OC
IV 1, p. 64)}

\begin{abstract}
 Simone Weil is one of the most prominent 20th century French philosophers. She is the sister of André Weil, the renowned  mathematician, the father of modern algebraic geometry and the initiator of the Bourbaki group. Simone and André Weil shared a love for  literature, mathematics, science and philosophy. My aim in this article is to convey, based on their writings and their correspondence, the idea that Pythagoreanism was a central element of their thought. I will put this into context, talking first about the life and work of each of them, showing how much they were linked by essential common ideas, even though their life paths were very different, and how, ultimately, Pythagorean mathematics and philosophy became naturally part of their respective intellectual worlds. The article is the written version of a lecture I gave in October 2025, at the conference 
 ``The Life and Contribution of Pythagoras to Mathematics, Sciences, and Philosophy" that  took place on October 3-4, 2025 at the Cyprus University of Technology in Limassol.

\bigskip

\bigskip


\end{abstract}
 \bigskip

 \section{Introduction}
 
  This article\footnote{This article is associated with a lecture I gave at the conference ``The Life and Contribution of Pythagoras to Mathematics, Sciences, and Philosophy" which took place on October 3-4, 2025 at the Cyprus University of Technology in Limassol.  I would like to thank Gregory Makrides, Sotos Voskarides and Christoforos Neophytidis for the invitation and for their care during the conference.} is centered on the role played by Pythagorean mathematics and mysticism in Simone Weil's thinking and writings. To put this in context, I will also talk about the life and work of Simone Weil and about her relation with her brother, the mathematician André Weil.
I start here by a few words on Pythagoreanism.

Among a multitude of quotations from modern authors on the influence of the Pythagoreans, let me quote Erwin Schr\"odinger (1887-1961), from his book \emph{Nature and the Greeks} \cite{Schrodinger}: ``The Pythagorean school is the prototype of a school of thinkers with
strongly scientific orientation and at the same time
with a well-marked bias, bordering on religious prejudice, towards reducing the edifice of nature to pure
reason". Let me also quote the mathematicians and historian of mathematics 
Michel Chasles (1793-1880). In his treatise \emph{Aperçu historique sur l'origine et le développement des
méthodes en géométrie, particulièrement de celles qui se rapportent à
la géométrie moderne} (Historical overview on the origin and development of the methods in geometry, in particular those relative to modern geometry) \cite{Chasles}, the name Pythagoras is amply mentioned. Chasles writes, on p. 4: ``It was primarily
to Pythagoras, who incorporated Geometry into his philosophy, and to his
disciples, that this science owed its first discoveries."\footnote{In this article, all the translations from the French are ours, except when an official translation is indicated in the bibliography.} 

The two quotes have a philosophical taste, and in the pages that follow, I will discuss indeed philosophy and spirituality, together with mathematics.
I will try to explain that with Simone Weil, like with the Pythagoreans, notions of harmony, proportion, beauty, rigor, and unity of opposites are no longer specific to mathematics, but govern the spiritual world. The Logos of the Pythagoreans acquires a mediating function, that of the second person of the Trinity.

From a purely mathematical point of view, let me recall that it is generally admitted that Book II of Euclid's \emph{Elements}, which several mathematicians, including André Weil whom we shall quote extensively, regard as a book on geometric algebra, is of Pythagorean origin. Concerning this matter, I refer the reader to the recent article by Negrepontis and Farmaki \cite{NF} in which the authors propose a restored Pythagorean version of that book.

The plan of the rest of this article is the following:

Section \ref{S-Weil-biography} contains some elements of Simone Weil's biography.

In Section \ref{s:A-Weil}, we mention some of the major mathematical achievements of André Weil.

Section \ref{s:relation-frere-soeur} is concerned with the relation between Simone and André Weil.

 Section \ref{s:Correspondence} is a preview on some discussions on mathematics and Pythagoreanism in the correspondence between the two Weils.

 In Section \ref{s:commentaries}, we review a few comments of Simone Weil on some Pythagorean texts, in which she expresses her belief that Pythagorean philosophical doctrine is close to that of Christianity. In fact, she considered the Pythagoreans as the spiritual ancestors of the Christians.

The final section, \S \ref{s:Bourbaki}, is concerned with Pythagorean mathematics as referred to in Bourbaki's \emph{Elements of history of Mathematics}. The relation between Simone Weil, André Weil, and Bourbaki is explained in the previous sections.

 \section{Simone Weil}\label{S-Weil-biography}
 
 There are a number of biographical writings on Simone Weil. The best biography, by far, remains, in my opinion, the one written by her friend Simone Pétrement \cite{Petrement}, to whom Simone Weil's parents entrusted the task of writing a biography of their daughter after her death, providing her with many documents they had at their disposal: notes, letters, etc. Other genuine information on Simone Weil's intellectual and spiritual development is found in her published correspondence, especially with her parents and her brother, as well as the latter's writings about his sister, particularly his memoirs published in English under the title \emph{Apprenticeship of a Mathematician}. There is also Simone Weil's correspondence with various people, in particular  the volume \emph{Attente de Dieu} (Waiting for God), a collection of letters she sent to  Fr. Joseph-Marie Perrin, who was her spiritual mentor, in which she writes at length about her philosophical position with respect to the Catholic Church. This collection includes a long    
 letter to Perrin that was later called \emph{Spiritual Autobiography} and of which we shall quote a passage below.
 Let me also mention the book \emph{Mon dialogue avec Simone Weil} (My dialogue with Simone Weil) written by the same Perrin, published 39 years after Simone Weil's death \cite{Perrin}.
 
 In the rest of this section, I will say a few words about Simone Weil which may give a brief overview of her intellectual and spiritual world, before getting to the heart of the matter.

  Simone Weil was born in Paris on April 3, 1909 to a father who was a doctor of Alsatian origin and a mother of Russian descent. In addition to her parents, she grew up under the loving and attentive gaze of her brother, André Weil, who was three years older than her and who became later the well-known mathematician. He taught her to read at a very young age. He writes, in his memoirs \cite[p. 22]{Weil-Souvenirs}:
``I devised a plan
to surprise my father for his birthday, April 7: my sister was to read the
newspaper to him".
André was eight years old, and Simone five. With her brother, she also learned Greek and developed a love of science. We shall say more about her brother later, especially in a section dedicated to him.

In 1928, Simone Weil was admitted to the \'Ecole Normale Supérieure, in the philosophy section, where she was the only woman in her class.\footnote{\label{s:ENS} For those unfamiliar with the French education system, let me recall that the \'Ecole Normale Supérieure, in Rue d'Ulm, in Paris, also known as ENS Ulm, founded in 1764, is the most difficult university institution to get into in France and probably in the world. Preparing for the ENS entrance exam usually takes two or three years after finishing high school, and studies there generally last four years. Very few students admitted every year. Until 2005, the school did not award any degree; graduates left with only the title: Ancien \'Elève de l'\'Ecole Normale Supérieure (Former Student of the \'Ecole Normale Supérieure). In 1928, Simone Weil entered the philosophy section of the ENS, at age 19. Her brother had entered, the mathematics section, in 1922, at age 16, which is very exceptional.} At the same time, she was enrolled at the Sorbonne, where she obtained in 1930 a ``diplôme d'\'Etudes Supérieures", with a dissertation titled:   \emph{Science et Perception dans Descartes} (Science and Perception in Descartes). In this work, she analyzes the question of whether science could bring equality or a certain degree of freedom to humankind, or, if, on the contrary, it leads humanity into a new form of slavery. These are themes that she was inclined to keep throughout her life. 
 To her dissertation, she had attached the epigraph   Ἀεὶ ὁ θεὸς γεωμετρεῖ.

During her years of study, 
Simone Weil was spending a substantial amount of her time and energy in political actions through distributing petitions, gathering money for workers' union strike and unemployment funds, and social needs. Even though she was firmly aligned with the left, she did not join any political party; the only label that could be applied to her is probably that of being an anarchist.

 The period of Simone Weil's university studies was also a time when Jean-Paul Sartre and 
Simone de Beauvoir were students. We know that she met them both, but had very little contact with them.
 Simone de Beauvoir, who was one year older than Simone Weil, relates in her \emph{Memoirs of a Dutiful Daughter} \cite{Beauvoir}
 that she sometimes caught sight of the latter in the courtyard of the Sorbonne with the newspaper l'Humanité in the pocket of her jacket.\footnote{The newspaper l'Humanité, founded by Jean Jaurès in 1906, was a left-wing daily from the outset.  In 1994, it became the official organ of the French Communist Party.} She writes: ``She intrigued me because of her reputation of intelligence and her bizarre attire [\ldots] A great famine had just devastated China, and I was told that upon hearing this news, she had burst into tears: these tears commanded my respect even more than her philosophical gifts. I envied a heart capable of beating for the entire universe".
   A little further, she writes:
 ``One day I managed to tame her. I don't remember how the conversation started, but she declared in a sharp tone that only one thing mattered in the world today: the Revolution, which would give everyone enough to eat. I replied just as emphatically that the problem was not to make people happy, but to find meaning in their existence. She looked me up and down and said, `It's obvious you've never been hungry.' That was the end of our relationship. I understood that she had labeled me a `spiritualist petit bourgeois,' and it irritated me." This was the last known interaction between the two would-be philosophers.

Simone Weil obtained her agrégation\footnote{The agrégation, in France, is a difficult national competitive exam, in principle for becoming a high school teacher, with  higher salary and fewer teaching hours than teachers without agrégation. Practically all students at the \'Ecole Normale Supérieure take the agrégation. Simone Weil came seventh in the agrégation in philosophy at age 22. Her brother André had  taken the agrégation in 1925 (in mathematics) at age 19, and was ranked first. Simone de Beauvoir came second in the philosophy agrégation in 1929, the same year as Jean-Paul Sartre, who came first. It was his second attempt; he had failed the first time. After getting the Agrégation, one is automatically assigned to a high school, unless the person wants to continue with doctoral studies and teach later at the university level.  Sartre and Beauvoir, like Simone Weil, became high school philosophy teachers, after obtaining the agrégation.} in 1931. She was assigned to teach philosophy in various high schools in provincial towns: first at the Lycée of  Puy-en-Velay, then in Auxerre, Saint Quentin and Bourges. The summer before she began her career as a philosophy teacher, she had taken a trip to Germany, from which she returned frightened by what she had seen of the rise of Nazism and the crushing of the labor movement. It was the year before Hitler came to power. She published several articles on this subject in philosophical journals and also in trade union brochures. Her articles from this period are collected in the volume \emph{\'Ecrits sur l’Allemagne (1932-1933)} (Writings on Germany (1932-1933)) \cite{Weil-Allemagne}.

During her period as a teacher, convinced that access to knowledge should not be a privilege reserved for the wealthy classes but should also be available to the unemployed and workers, Simone Weil used to organize courses and lectures for the latter. Very early on, she began to donate a large part of her salary to unemployment funds.  She came into conflict with the administration of her school, which disapproved her social activism and her commitment to defending the poor and the working class.  In 1934, she  took a leave from her job as a teacher with the intention of becoming a factory worker. Her aim was to physically experience the social oppression, the fatigue of assembly line work, and the constant anxiety of losing one's job that laborers of that era were subjected to. 
She worked at several factories, always asking to be assigned to the most strenuous jobs.  Her health was gradually deteriorating, aggravated by voluntary hunger. During this time, she was publishing articles on social problems, like the status of women in factories, as well as political pamphlets on various subjects, including one on the contradictions in Marxism, even though she was very familiar with Marx' work and admired it. At the same time, she was publishing essays on Ancient Greek and Roman history, on the Iliad and on classical Greek theatre, some of them addressed to her colleagues in the factory. She had the conviction that the latter are capable of understanding these tragedies, as much as the educated class.  

In 1942, Simone Weil fled France with her parents, with the aim of protecting them, after the Vichy laws were introduced. Her brother André Weil was already in America, where he was offered a job at Lehigh University in Bethlehem (Pennsylvania). After her parents were settled in the United States, she tried to return to France to join the resistance, but she wasn't allowed to because of her poor health. She ended up in England, where she joined the France Combattante.\footnote{The France Combattante was a union of France organisations who were fighting against the Vichy regime and the German occupation of France. The movement  was recognized by Great Britain, where Simone Weil was stationed, waiting to turn back to France.} She was assigned there to an office.
She died a few months later, on August 30, 1943, in a hospital in the suburbs of London. She was being treated there for tuberculosis, but the cause of her death was malnutrition.

 From a certain point on, Simone Weil became totally immersed in Christian mysticism. From a letter to Father Perrin \cite[p. 121]{Weil-Attente}: ``Peasants and workers possess that incomparably sweet closeness to God which lies at the bottom of poverty, lack of social consideration and long, slow suffering." In order not to overstep the bounds of this paper, I will touch only briefly on Simone Weil's mystical side, even though it is an essential part of her thinking and is closely linked to her deep interest in ancient Greece, Pythagoreanism, and mathematics. I will mention her religious thinking in passing, in relation with Pythagorean mathematics, in Section \ref{s:commentaries}. Readers interested in Simone Weil's Christian mysticism and its relationship with mathematics can also consult our article 
 \cite{Papadopoulos-Weil-Pythagoreans}. Simone Weil's mysticisim and her ideas on mathematics are also discussed in the article \cite{Lafforgue} by Laurent Lafforgue.
 
In her biography, Simone Pétrement recounts that Simone Weil's mother, responding to a person who was talking to her about her daughter's glory, said: ``Ah! I would have so much preferred that she had been happy!''\cite[Vol. 1, p. 20]{Petrement}.

 It is time to say a few words about André Weil, Simone's brother.

\section{André Weil} \label{s:A-Weil}
André Weil is one of the most important mathematicians of the twentieth century. We have already mentioned his relationship with and influence on his sister. As a child, the least we can say is that he was intellectually precocious. Simone Pétrement, the official Simone Weil biographer, also talks about the brother, based on their mother's memoirs and other family letters. She reports that in 1914 (André was 8 years old), he came across a geometry book by \'Emile Borel, began to study it, and was able to solve difficult exercises in the book. With this, the suspicion of his mathematical talent became an evidence \cite[vol. 1, p. 26]{Petrement}.

André Weil  is the main founder of modern algebraic geometry, and he made substantial advances in several other domains.
 In the lines that follow, I will mention some of the fields in which he made important contributions to mathematics.
  
The name André Weil is attached to a number of objects and theories that had a profound and durable impact on mathematics. Among these, there are the Borel--Weil--Bott theorem in the representation theory of Lie groups,  the Bergman--Weil formula for the representation of holomorphic functions of several variables---an integral formula that generalizes to several variables  the fundamental Cauchy integral for holomorphic functions of one variable---, the Shafarevich--Weil theorem in Galois theory, the de Rham--Weil theorem  in algebraic topology,  Chern--Weil theory  and in particular the Chern--Weil homomorphism in the theory of characteristic classes,   the
  Weil--Petersson K\"ahler metric on Teichm\"uller space, 
 the introduction of the metaplectic group and a representation of this group, called the Weil representation, the Mordell--Weil theorem in algebraic number theory saying that the group of points on an abelian variety that are rational over a given number field is always of finite type,
  and the Taniyama--Shimura--Weil conjecture connecting topology and number theory, which is settled now and which is called the Modularity theorem, of which the Fermat last theorem is a special case. This is only part of the list of mathematical concepts called after André Weil.

  I will say no more about André Weil's work, except to point out that his favourite field was arithmetic, and most of his books and articles are on this topic. Let me mention, without entering in the details, the titles of  his doctoral thesis, defended in 1928,  \emph{Arithmetic on Algebraic Curves} \cite{Weil-1928} and that of his his 1935 booklet, \emph{Arithmetic and Geometry on Algebraic Varieties} \cite{Weil-1935}. This book was followed by a note published in Russia, titled \emph{Arithmetic of algebraic varieties} \cite{Weil-1937}. Let me also mention that
in  1950, at the height of his career, André Weil gave a plenary lecture at the ICM in Cambridge, Massachusetts, bearing  almost the same title as his article in Russian, \emph{Number theory and algebraic varieties} \cite{Weil-1950}.
    I am emphasizing arithmetic here because it was the favorite field of the Pythagoreans, whose motto was that ``everything is number" (meaning an integer).

During his stay at the University of Strasbourg (1932-1939), André Weil lectured on arithmetic. It is interesting to read what he says about this, in the \emph{Apprenticeship}: 
``I taught a course on
algebraic number theory. I believe this was the first number theory course
to be offered in a French university since the beginning of the century. In
the university catalog, it was to be listed as  `Arithmetic.' To the dean, this
title smacked of primary school; it did not correspond to his notion of the
university's honor, so we changed the course title to `Number Theory' or
perhaps `Higher Arithmetic,' and he was happy" \cite[p. 111]{Weil-Souvenirs}.

  It is also important for our purposes here to remember that André Weil
 is the main founder of the Bourbaki group. Cartier writes in \cite{Cartier}: ``Bourbaki was undoubtedly Weil's creation, even if there was no shortage of strong personalities." Bourbaki's \emph{Elements of mathematics}  \cite{Bourbaki-Ensembles} is full of references to the Greeks, in particular the Pythagoreans. It is very probable that a large number of these references are due to Weil. Among the Bourbaki members he was the one with the most profound knowledge of history.
  I will quote a few passages on this matter, in \S \ref{s:Bourbaki} of the present article.  
      
       It is useful to keep in mind that Bourbaki was founded  with the goal of writing a treatise  that would be an analogue of Euclid's \emph{Elements}, first for analysis, and then the project was extended for all mathematics.\footnote{The first members chose as a collective name for the group that of a French general of Greek descent, Constantin Denis Bourbaki (Κωνσταντῖνος Διονύσιος Βούρβαχης). The latter, along with other French officers, left for Greece with the permission of the French government to fight the Ottoman occupation after the Greek revolution of 1825 broke on. Bourbaki arrived in Greece around the end of November 1826, leading a Philhellene Corps. Bourbaki was devoted to Karaiskakis. Makrygiannis mentions him in his Memoirs \cite{Makriyannis}. He died in Greece, on January 27 1827, at age 40.
His name is mentioned on an engraving honoring the Philhellenes on display at the Benaki Museum in Athens (on a staircase).
Several members of his family, including his eldest son, Denis, were also fervent Philhellenes. The first name of Bourbaki (the mathematician), Nicolas, was chosen by his godmother, Eveline, the future wife of André Weil. André Weil writes: ``My future
wife Eveline, who was present at this discussion, became Bourbaki's
godmother and baptised him Nicolas." \cite[p. 81]{Weil-Souvenirs}.}   Bourbaki shares with the Pythagoreans, in addition to his dedication to mathematics, a cult of mystery: 
  the names of the members were kept secret until after their retirement from the group. 
  Let us also recall that like the Pythagoreans of southern Italy at their peak, Bourbaki wielded at some point an undeniable occult power in the French academic landscape.

Simone Weil  participated in the first Bourbaki congresses; see the photos in Figures \ref{fig:Bourbaki-1937} and  \ref{fig:Bourbaki-1938}  below. We shall insist on the fact that she understood and retained several issues that were debated during these meetings, and that she later discussed them with her brother.

 \section{On the relation between the brother and sister}\label{s:relation-frere-soeur}

In an interview with André Weil, concerning his relation with his younger sister, and to the question: ``Was there a strong intimacy between the two of you?", he responded: ``A very big one. As a child, my sister would spend her time imitating me" \cite[p. 10]{Weil1990}. Simone Pétrement quotes in her biography \cite[Vol. 1, p. 24]{Petrement}, from a letter written by Mme Weil to a friend of her about her daughter: 
``Simone has developed in an incredible way. She follows André everywhere, takes an interest in everything he does, and, like him, now finds that the days are too short. They have an excellent influence on each other: he protects her, helps her climb difficult passages, often gives in to her, and she, who is with him from morning to night, becomes livelier, happier, and more enterprising."
André Weil writes in the Foreword to his book \emph{The Apprenticeship of a Mathematician} \cite[p. 11]{Weil-Souvenirs}:    ``As children, Simone and I were inseparable; but I was always the big brother and she the little sister. Later on we saw each other only rarely, speaking to one another most often in a humorous vein; she was naturally bright and full of mirth, as those who knew her have attested, and she retained her sense of humor even when the world had added a layer of inexorable sadness.
In truth we had few
serious conversations. But if the joys and sorrows of her adolescence were
never known to me at all, if her behavior later on often struck me (and
probably for good cause) as flying in the face of common sense, still we
remained always close enough to one another so that nothing about her
really came as a surprise to me---with the sole exception of her death. This I did not expect, for I confess that I had thought her indestructible. It was not until quite late that I came to understand that her life had unfolded according to its own laws, and thus also did it end. I was little more than a distant observer of her trajectory."

\bigskip

 \begin{figure}
\centering
\includegraphics[width=15cm]{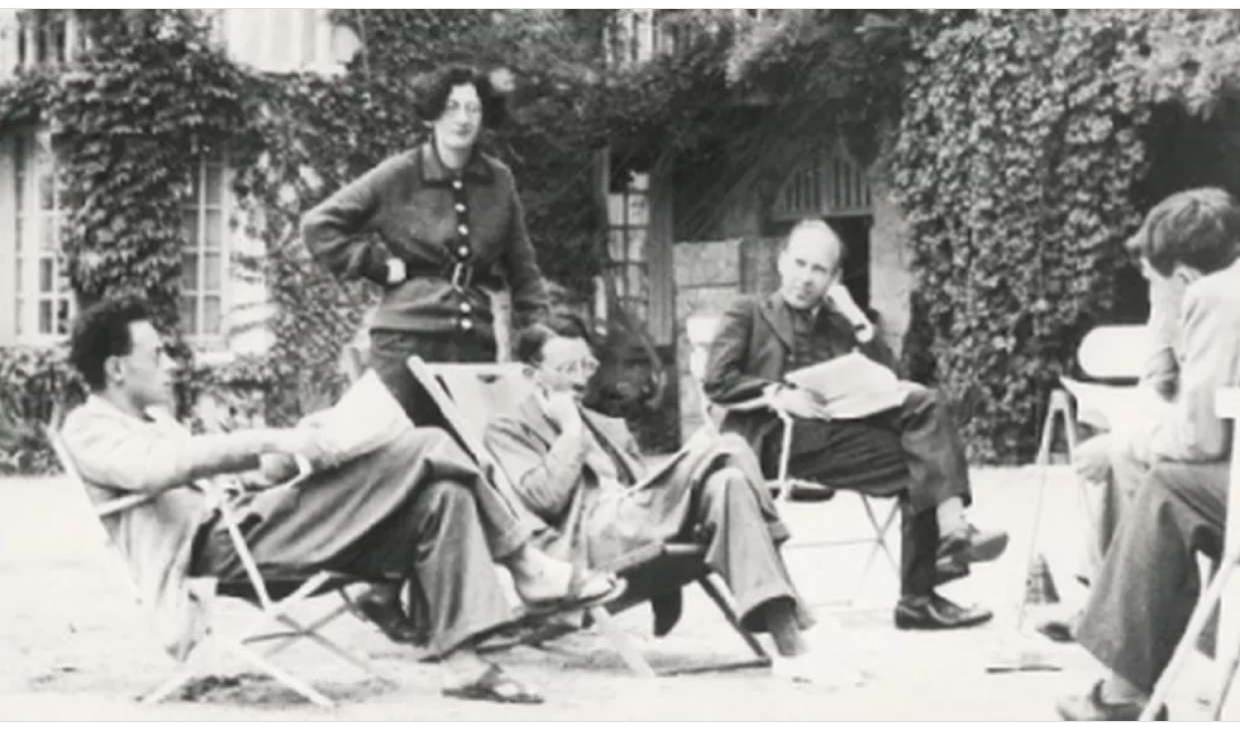}
\caption{A group photo taken at the Bourbaki congress at Chan\c cay, 1937. From left to right, André Weil, Henri Cartan, Szolem Mandelbrojt, and  Simone Weil standing. André Weil writes in \cite[p. 118]{Weil-Souvenirs}: ``In the meantime, in September of
1937, Bourbaki was again to meet in Chan\c cay. My sister, still intent on
improving her knowledge of mathematics, attended our congress". Chan\c cay is a village  situated in the French department of  Indre-et-Loire.}
 \label{fig:Bourbaki-1937}
\end{figure}

 \begin{figure}
\centering
\includegraphics[width=15cm]{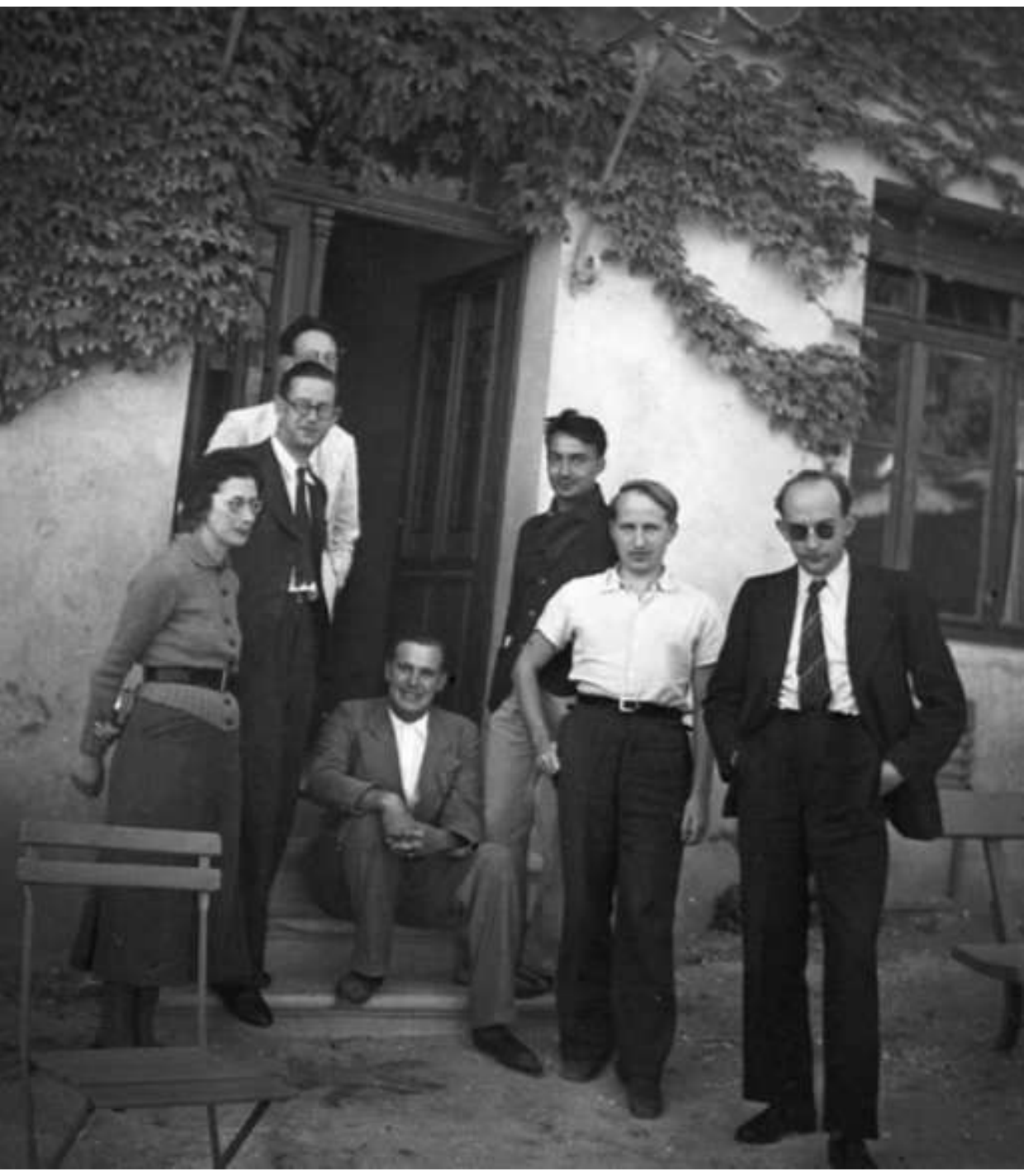}
\caption{Bourbaki congress at Dieulefit,  Septembre 1938. From left to right: Simone Weil, Charles Pisot,  André Weil, Jean Dieudonné (sitting), Claude Chabauty, Charles Ehresmann and Jean Delsarte.  Dieulefit is a commune in the   department of Drôme, in the region called Auvergne-Rhône-Alpes, Southeast of France.  The photo was taken for the third anniversary of Bourbaki.}
 \label{fig:Bourbaki-1938}
\end{figure}
  
 \bigskip

In a famous passage of a letter that Simone Weil sent to J.-M. Perrin,  
 published under the title
 \emph{Autobiographie spirituelle} (Spiritual autobiography) in the volume \emph{Attente de Dieu} \cite{Weil-Attente}, she writes about her brother: 
 ``At age fourteen, I fell into one of those bottomless despairs of adolescence, and seriously thought about dying, because of the mediocrity of my natural faculties. 
The extraordinary gifts of my brother, who had a childhood and youth comparable to Pascal's, forced me to be aware of this. I had no regrets about external successes, but about not being able to hope for any access to that transcendental realm where only truly great men can enter, and where truth dwells. I would rather die than live without it. 
After months of inner darkness, I suddenly and forever had the certainty that any human being, even if his natural faculties are almost nil, can enter that realm of truth reserved for genius, if only he desires the truth and makes a perpetual effort of attention to reach it. Later, when headaches weighed heavily on what little faculties I possess, causing a paralysis that I quickly assumed would probably be permanent, that same certainty made me persevere for ten years in efforts of attention that were supported by almost no hope of results."

In the article \cite{Papadopoulos-Tanabe}, we discuss the life and thinking of Simone and André Weil, emphasizing their attachment to India and Indian mysticism, another important aspect of the bond between Simone and André Weil.

  We now turn to Simone Weil's correspondence with her brother

 \section{Correspondence on mathematics and the Pythagoreans}\label{s:Correspondence}
   Several letters exchanged between André and Simone Weil reveal the latter's desire to understand her brother's work, and in particular, his results in number theory. Conversely, André Weil was interested in his sister's comments on his writings. 
In a letter sent from Paris, dated February 17, 1934,
 after informing her of arithmetical problems on transcendence, he writes: ``I will wait patiently for the development of your ideas on mathematics"  \cite[p. 529]{SWeil-Corresp}.

When it came to questions about history of mathematics, and in particular on the relationships between Greek, Babylonian, Hindu and other mathematical schools, Simone Weil had her own ideas; she was able to establish logical connections, and her brother was pleased and interested in hearing them. He enjoyed their correspondence, particularly on this topic.  In the letter to his wife written from prison,\footnote{Having refused to fulfill his military service, André Weil was first sentenced to five years in prison. Eventually, he was imprisoned in France from February to May 1940, first in Le Havre and then in Rouen. This is another story, which I discuss elsewhere, see \cite{Papadopoulos-Tanabe} and \cite{Papadopoulos-Weil-Pythagoreans}.}  that we quoted above \cite[p. 150]{Weil-Souvenirs}, he writes:  
  ``[\ldots] And then there is my correspondence with my sister; nothing could
be easier than this at the moment, since it concerns only the most abstract
subjects. Greek mathematics, for example: she has some ideas about that." Several times, his sister would ask him for mathematical references and writings by historians of mathematics that would allow her to test her ideas. A few times, she enquired about  Neugebauer,\footnote{Otto Neugebauer (1899-1990) was a mathematician and historian of mathematics and astronomy. He was a specialist, among other things, of mathematics of ancient Mesopotamia. His name is mentioned several times in the correspondence between André and Simone Weil, and in André Weil's book \emph{Number theory, An approach through history. From Hammurapi to Legendre} \cite{Hammurapi}.} the major specialist of Babylonian mathematics, asking for references to his works. The possible transmission of ideas on number theory from Mesopotamia to Greece that she alluded to is accounted for in the chapter on ``Protohistory" of Weil's book \emph{Number theory, An approach through history. From Hammurapi to Legendre} \cite{Hammurapi}.
In a letter dated February 1940,  sent to her brother in prison and published in her \emph{Works}, t. VII vol. 1, Simone Weil, after describing  to him some passages she read in Neugebauer's book \cite{Neugebauer-1934} on Babylonian mathematics, she writes, about the latter \cite[p. 437]{SWeil-Corresp}:
\begin{quote}
\small
Strange people, those Babylonians. Personally, I am not very fond of their abstract thinking---the Sumerians must have been much more pleasant. After all, they invented all the Mesopotamian myths, and myths are much more interesting than algebra. But you must be a direct descendant of the Babylonians. As for me, I believe that God, according to Pythagorean teachings, is a perpetual geometer---but not an algebraist. In any case, I was pleased to read in your last letter that you deny belonging to the abstract school.
   \end{quote}
    André Weil was interested in his sister's ideas. He wrote to her: “I find your ideas on the role of proportion in the history of Greek thought very appealing and, at the same time, quite plausible, [\ldots]” (letter dated March 28, 1940 \cite[p. 556]{SWeil-Corresp}) .
    
    I would like to highlight two points in Simone Weil's passage I just quoted. The first is her reference to God and to Pythagorean teaching. This brings us to the fact that Simone Weil, in her philosophical writings, repeatedly stated that she considered the ancient Greeks, starting with the Pythagoreans, to be the true precursors of Christians. In claiming this, she departs from the official position of the Church, which considers that it was the Jewish religion that heralded the coming of Christ. In Section  \ref{s:commentaries}, we shall touch on this question, but to discuss it in depth would take us far from the scope of this article. We discuss this in the article \cite{Papadopoulos-Weil-Pythagoreans}.
        
           The second point which I would like to highlight is the sentence: ``I was pleased to read that you deny belonging to the abstract school". This fact goes against a certain reputation of Bourbaki, and of André Weil, of being too formalistic and too abstract. On the same subject, we can quote a passage from a letter written in prison by André Weil to his wife that he mentions in the \emph{Apprenticeship}, dated April 22, 1940 \cite[p. 150]{Weil-Souvenirs}: 
``My mathematical fevers have abated; my conscience
tells me that, before I can go any further, it is incumbent upon me to work
out the details of my proofs, something I find so deadly dull that, even
though I spend several hours on it every day, I am hardly getting anywhere."

\begin{figure}
\centering
\includegraphics[width=7cm]{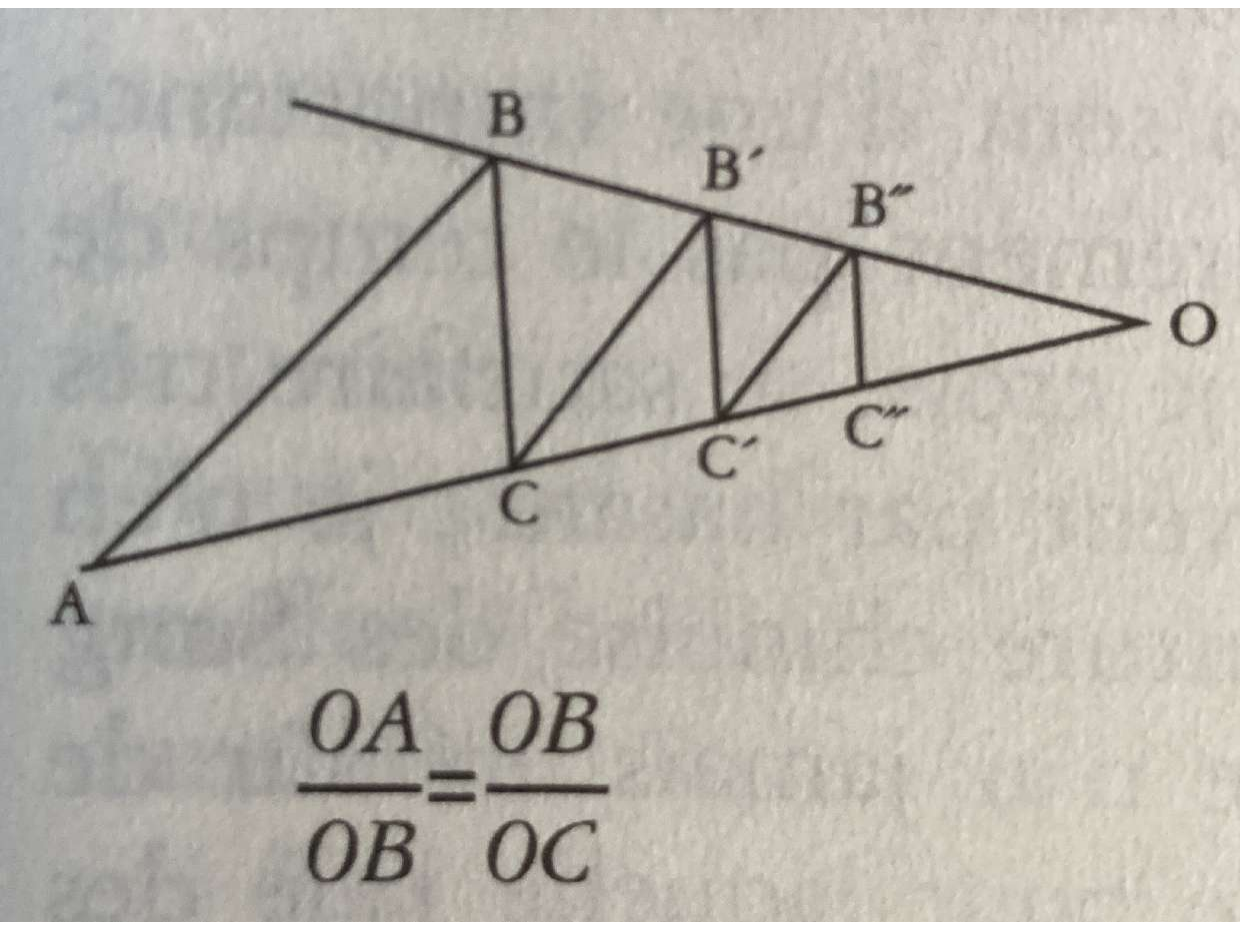}
\caption{From a letter sent by Andr\'e Weil to his sister, dated March 28, 1940, commenting on the Pythagorean proof of the incommensurability of the diameter to the side of a square, \cite[p. 554-555]{SWeil-Corresp}.}
 \label{fig:Weil}

\end{figure}

   In another letter, sent by Andr\'e Weil to his sister, dated March 28, 1940, he writes \cite[p. 554-555]{SWeil-Corresp}: ``\emph{How} do you think [the Greeks] discovered the immeasurable? I remember that there is no sure data on this, and that I saw somewhere an ingenious suggestion, with an infinite descent\footnote{Infinite descent  is a method based on the fact that a strictly decreasing sequence of natural numbers is necessarily finite.  
  André Weil was fond of this method because it forms the basis of the proof of Mordell's theorem, a special case of what became later the Mordell--Weil theorem.  Besides, this method was a crucial argument in André Weil's doctoral thesis (1928), in which he introduced a notion of height that allows to bound the size of $K$-rational points in the group of rational points in an abelian variety over a number field $K$.
 We can read in this thesis \cite{Weil-1928}:
 ``This method,
applied systematically for the first time by Fermat, who gave it its name,
consists, as we know, of providing a process by which, from any solution
to an equation under study, another can be deduced, and of showing that the iteration of this process cannot be continued indefinitely. This is how, for example,
Fermat proved the impossibility of $y^2=x^4-z^4$ in integers, by showing
that from any solution, another one can be deduced in smaller integers." Weil's book on \emph{Number theory: An approach through history} contains an appendix titled \emph{The Descent and Mordell's Theorem}  \cite[p. 140-149]{Hammurapi}.} 
 where the thing could be \emph{seen} on a square and its diagonal (numerically this gives the following: We have $\frac{\sqrt{2}+1}{1}= \frac{1}{\sqrt{2}-1}$; but if $\frac{a}{b}=\frac{c}{d}$ and $a,b,c$ have a common measure,  
then the same holds for all the lengths obtained by continuing the progression downwards, which is absurd, since these lengths $OB', OC', OB'', OC''$, etc. become as small as we want)" (see Figure \ref{fig:Weil} here). In fact, it is admitted that the Pythagoreans discovered the irrationality of $\sqrt{2}$ by showing that the anthyphairetic expansion of the diagonal to the side of a square is infinite, an operation which is vaguely related to, but much simpler, than the method of infinite descent (both methods are based on the property that of the set of natural numbers is well-ordered). See the recent article \cite{NP2026}. A little bit further in the letter, André Weil writes: ``It is probable that, even for $\sqrt{2}$, people started in this way by something equivalent to the decomposition in prime factors." He concludes:
``In any case, I know of nothing that contradicts the idea (which I no longer remember where I got it from, but which is not my own)\footnote{One must recall that in prison, André Weil did not have his documents or a mathematical library available.} that it was the discovery of irrational numbers that led to the first axiomatic research and the first \emph{proofs}\footnote{The emphasis is André Weil's.} in the modern sense."

 \bigskip

 To André Weil, who was in prison, his sister sent another letter, in February 1940 (with no precise day), in which she writes \cite[p. 436]{SWeil-Corresp}: 
  
 \begin{quote} 
 
Who knows, maybe you will discover some fascinating things? But there is another distraction now that you have free time. I don't know if I mentioned it in the letter I wrote to you in Le Havre, which you must have by now, but never mind. It would be a way of making ordinary people (like me, for example) understand what your current research is all about. I am sure it would be very good practice for you. What do you have to lose? You are not going to waste your time, because you have time to spare. It is all very well to mock those who, like a friend of mine from my days at Rue d'Ulm\footnote{\'Ecole Normale Supérieure, see Footnote \ref{s:ENS}).}, philosophize about mathematics without knowing anything about it, but perhaps it is precisely  mathematicians who should try to do this work. [\ldots] What would it cost you to try? I would find it fascinating.

 \end{quote}

Simone Weil repeatedly urged her brother to write an overview, for non-specialists, of his current research and its significance. André Weil finally responded positively, sending her a long text in the form of a letter dated March 26, 1940, in which he explains his arithmetico-algebraic work (the term he uses). The reflections contained in this letter are of two kinds, he says: some concern the history of number theory, and others the role of analogy in mathematical discovery. It is from this second part that he thought his sister would benefit most.

This letter, which has become a reference text, is published in Simone Weil's correspondence \cite[p. 535-552]{SWeil-Corresp}, and also in the Complete Works of André Weil, see \cite{Weil-Works}, text No. [1940a].

Meanwhile, Simone Weil continued to constantly reflect on mathematics and on history of mathematics.  
In a letter to her brother (who was already in the United States), written in April 1941 (two years before her death),  she writes \cite[p. 489]{SWeil-Corresp}: 
``Can you imagine that recently, while trying to rediscover Archimedes' ‘mechanical method’ for squaring the parabola, I found another method similar to his, which he could very well have used as well!" and she explains her method. In the same letter, she writes : 
``Can you let me know if Neugebauer has published any new works on ancient mathematics or astronomy since the one I had in my hands (which I believe dates from 1934)?" 
Regarding Indian mathematics, she writes to him in September 1941 
 \cite[p. 498]{SWeil-Corresp}:  ``Reading the Upanishads, I am struck more and more by the similarity between Indian and Greek thought, to the point that I am almost inclined to believe in an Indo-European philosophy that would be the common source; for they are similar in certain respects to the point of identity, without bearing any trace of mutual influence."

 \section{On Simone Weil's commentaries on Pythagorean fragments}\label{s:commentaries}

We already alluded to the fact that Simone Weil considered that the  mystical foundations of Christian theology can be found in the writings of the Pythagoreans. 
In her collection of writings that were published under the title \emph{La source grecque} (The Greek source), we find translations and commentaries of several passages of Greek texts, including fragments from the Pythagorean Philolaos.

She translates the fragment 

ἔστι ἁρμονία πολυμιγέων ἕνωσις καὶ δίχα φρονεόντων συμφρόνησις

 \cite[p. 155]{Weil-Source}\footnote{In Nicomachus of Gerasa's \emph{Introduction to arithmetics}, Diels, I, 410, fr. 10.}  as follows: ``Harmony is the unity of a mixture of many, and the unique thought of separate thinkers."\footnote{Her French translation: ``L’harmonie est l’unité d’un mélange de plusieurs, et la pensée unique de pensants séparés."}
She sees in this fragment a theological connotation. She writes:  ``It seems to me that the second part of the definition can only apply to a being in several persons," and she compares this fragment to another fragment of Philolaos (DK58B1), 

φιλίαν ἐναρμόνιον ἰσότητα.

which she translates as follows\footnote{Diels, 451, lines 12-13.} ``Friendship is an equality built on harmony",\footnote{``L’amitié est une égalité faite d’harmonie."}. She declares that  the two formulae combined would be a perfect starting point for a theologian who wants to speak about love in the Trinity.
A little further, she translates the fragment DK44B11: 

ἴδοις δὲ οὐ μόνον ἐν τοῖς δαιμονίοις καὶ θείοις πράγμασιν τὴν τοῦ ἀριθμοῦ φύσιν καὶ τὴν δύναμιν ἰσχύουσαν, ἀλλὰ καὶ ἐν τοῖς ἀνθρωπικοῖς ἔργους καὶ λόγοις πᾶσι παντᾶ καὶ κατὰ τὰς δημιουργίας τὰς τεχνικὰς πάσας καὶ κατὰ τὴν μουσικήν.

``We see that the essence and virtue of number do not reign \emph{only} among religious and divine things, \emph{but also} in all human actions and relations and in everything that relates to the technique of workmanship and to music.\footnote{``On voit que l’essence et la vertu du nombre ne règne pas \emph{seulement} parmi les choses religieuses et divines, \emph{mais aussi} dans toutes les actions et relations humaines et dans tout ce qui a rapport avec la technique des métiers et avec la musique." (The emphasis is Simone Weil's.)}" Commenting these fragment, she writes that they show that the true meaning of early Greek mathematics,  which is the foundation of our science, was religious. She declares that this is confirmed by the following passage from Plato's \emph{Epinomis}, 990d: 

\begin{quote}\small

[\ldots that] which is ridiculously called the measurement of the earth,\footnote{Plato talks here about geometry.} and which is nothing but the assimilation of numbers not naturally similar\footnote{Non commensurable.} to one another, made manifest by the purpose of plane figures. It is evident to anybody who can understand that this marvel is of divine, not human, origin."\footnote{``[\ldots] ce qu’on appelle ridiculement mesure de la terre, et qui n’est que l’assimilation des nombres non naturellement semblables entre eux, rendue manifeste par la destination des figures planes. Il est évident pour quiconque peut comprendre que cette merveille est d’origine non pas humaine, mais divine."}
\end{quote}
She comments on this passage as follows: 
 
 ``This means that geometry is the science of finding proportional means through the incommensurable proposition (i.e., $\sqrt{m\times n}$), and that it proceeds from a supernatural revelation. This, when considered alongside the passages where Plato describes the mediation between God and man through the image of the proportional mean ($a/b=b/c)$ and the numerous passages in the Gospels where Christ uses the same image for his own mediating function (``As my Father has sent me, so I send you, etc."), makes the invention of geometry in Greece, in the 5th century BC, appear as a prophecy in the strictest sense."
 Lafforgue, in his paper \cite{Lafforgue}, recalls that in St Paul's letter to the Hebrews, the name ``Mediator" is explicitly given to Christ (He 9,15 ; 12,24).
Simone Weil continues: 
``Perhaps the least inaccurate translation of the beginning of the Gospel of St. John would be:
\begin{quote}\small
In the beginning was the Mediation, and the Mediation was by the side of God, and the Mediation was God. It was in the origin by the side of God. All things came into being through it; without it nothing came into being. In it was life, and life was the light of mankind; and light shone in darkness, and darkness did not overcome it. It was the true light who came into the world, that illuminates every man. And it found itself in the world and the world was born through it, yet the world did not recognize it."
\end{quote}

She writes, a little further:
   ``[\ldots]  Thus, divine Mediation, through analogical descent, penetrates everything. It unites God to God, God to the world, the world to itself; it constitutes reality in all domains."
   
   We are back to the central figure of Mediation, that is, Jesus Christ. A mediation
allows us to connect different beings and different entities. Christ, for her, is the Mediation between humanity and the Father, and between the Father and the Holy Spirit. She adds that ``all this is expressed in the single term `logos',  as the name of the second Person of the Trinity." 
At several points, she declares that the term  λόγος, by which St. John,  in his gospel, designates the second Person of the Trinity, is synonymous to Mediation.

 \section{André Weil, history of mathematics, Bourbaki, formalisation and the Pythagoreans again}\label{s:Bourbaki}

André Weil's contribution to the history of mathematics (and in particular the history of Greek mathematics) is enormous.
He had strong ideas on who can write on history, pointing out the sometimes disastrous influence of those who write on the history of mathematics without being mathematicians. Regarding this matter, he writes in \cite{Weil-Letter}: ``when a discipline, intermediary in some sense between two
already existing ones (say A and B) becomes newly established, this often makes
room for the proliferation of parasites, equally ignorant of both A and B, who
seek to thrive by intimating to practitioners of A that they do not understand B,
and vice versa. We see this happening now, alas, in the history of mathematics.
Let us try to stop the disease before it proves fatal."   His book  \emph{Number theory. An approach through history. From Hammurapi to Legendre} 
\cite{Hammurapi}  is both  on mathematics and the history of mathematics. In it, we find the development of his idea that ancient number theory is an ancestor of modern algebraic geometry.

We have already talked about the pre-eminent place that André Weil occupied in the Bourbaki group.
The latter's volume \emph{Elements of the History of Mathematics} \cite{Bourbaki-H} is, broadly speaking, a compilation of the historical notes contained in the various volumes of the whole \emph{Elements of Mathematics}. We shall quote now a few passages from it. 
 Naturally, when discussing the origin of mathematics, Bourbaki gives the exclusive place to the Greeks. 
At the beginning of the first chapter of his \emph{Theory of sets} (the first book of the \emph{Elements of mathematics}), he writes:

\begin{quote}\small 
Since the Greeks, mathematics has been synonymous with proof; some even doubt that, outside of mathematics, there are proofs in the precise and rigorous sense that the Greeks gave to the word and that we intend to give it here. It is fair to say that this meaning has not changed, for what was a proof for Euclid is still a proof in our eyes; and in times when the concept threatened to be lost and mathematics was therefore in danger, it was among the Greeks that models were sought.
\end{quote}

In the  \emph{André Weil Fund} \cite{Andre-Weil-Fonds},
there is  a draft of the introduction to this book, in which we find a variation on the passage we just quoted, with some more details. It reads:

 \begin{quote}\small

 There is one mathematics, indivisible: this is the \emph{raison d'être} of this treatise, which aims to explain its elements in the light of a tradition spanning twenty-five centuries. However, since the Greeks, mathematics has meant proof; it is even doubtful whether proofs exist or could exist outside mathematics in the precise and rigorous sense that the Greeks gave to the word and that we intend to give it here. We can rightly say that this meaning has not changed, for what constituted a proof for Euclid is still a proof in our eyes; and in times when the concept threatened to be lost and mathematics was therefore in mortal danger, it was among the Greeks that its paradigms were sought.  
  \end{quote}
   The expression ``one and indivisible" has obvious theological connotations.
 
 Bourbaki considered the return to the axiomatic method to be precisely a return to the Greeks.  He writes in the \emph{Elements of the History of Mathematics} \cite[p. 10]{Bourbaki-H}:
 
  \begin{quote}\small

The essential originality of the Greeks consists precisely in a conscious effort to arrange mathematical proof in such a sequence that the passage from one link to the next leaves no room for doubt and commands universal assent.
That Greek mathematicians used, in the course of their research, just like modern mathematicians, heuristic rather than probative reasoning is what Archimedes' treatise on method, for example, would prove (if needed).
We will also note,
in the latter, allusions to results (found, but not proved)
 by earlier mathematicians. But, from the first
detailed texts known to us (and which date from the middle of the
5th century), the ideal ``canon" of a mathematical text is well established.
It found later its most complete realization among the great classics,
Euclid, Archimedes and Apollonius; the notion of proof, in
these authors, differs in no way from ours.
We have no text that allows us to follow the first
steps of this ``deductive method", which already appears close to perfection at the very moment we observe its existence.
We can only assume that it fits quite naturally into
the perpetual search for ``explanations" of the world, which characterizes
Greek thought and which is already so visible among the Ionian philosophers
of the 7th century; moreover,   tradition is unanimous in attributing the development
and refinement of the method to the Pythagorean school,
at a time between the end of the 6th and the middle of the 5th century.

\end{quote}
 
  The \emph{axiomatic method}, writes Bourbaki at the beginning of the \emph{Theory of Sets}, ``is, strictly speaking, nothing more than the art of writing texts whose formalization is easy to conceive." We note that the expression ``whose formalization is easy to conceive" is not identical to`` formal."
 It is also important to remember that Weil did make the difference between texts written in a formalized way and those written for research, and this is something that is not sufficiently emphasized by those who speak of him. After discovering a new result, he hated the step where he had to write everything down formally. We already quoted a passage from a letter to his wife in which he writes that he finds working on the details of his proof to be  ``deadly dull" \cite[p. 150]{Weil-Souvenirs}.

We have quoted in \S \ref{s:Correspondence} a discussion between André Weil and his sister on the  incommensurability of the side of a square with its diagonal (that is, the irrationality
of $\sqrt{2}$).   
In the chapter on \emph{Real Numbers} of the \emph{Elements of history of mathematics}  \cite[p. 185]{Bourbaki-H}, Bourbaki talks about this result, which he naturally attributes to
 the Pythagorean school. He writes:
 \begin{quote}\small
 Early Greek mathematics
was inseparably linked with speculations, partly scientific,
partly philosophical and mystical, on proportions, similarities
and ratios, in particular ``simple ratios" (expressible
by fractions with small numerators and denominators);
and it was one of the characteristic tendencies of the Pythagorean school
to claim to explain everything by the integer and the ratios
of integers. But it was the Pythagorean school, precisely, who discovered
the incommensurability of the side of a square with its diagonal (the irrationality
of $\sqrt{2}$): the first example, without doubt, of a proof
of impossibility in mathematics; the mere fact of asking such a
question implies the clear distinction between a ratio and its approximate
values, and is enough to indicate the immense gulf that separates Greek mathematicians
from their predecessors.
 \end{quote}
 
André Weil was convinced that the Greeks had relations with the Babylonians and were familiar with their mathematics, like his sister who had written to him many times on this matter (see the letter of February 1940 already quoted in \S \ref{s:Correspondence} above, \cite[p. 435]{SWeil-Corresp}).  But he also shared with his sister the conviction that the Greeks were not particularly interested in the kind of mathematics the Babylonians developed because what the Greek were looking for was not results but the rigor of the proof.
Their discovery of irrationality took on a religious character, because this existence was derived from a rigorous proof.  
Let us listen to what Bourbaki writes on this subject:    
 In the same  \emph{Elements of history of mathematics}, after talking about the fact that the Pythagoreans tried to relate everything to number, Bourbaki writes (p. 35):  
``[\ldots] Although the discovery of irrational numbers seemed to close this path forever,
the reaction it triggered in Greek mathematics
was a second attempt at synthesis, this time taking
geometry as its basis and incorporating, among other things, the methods of
solving algebraic equations inherited from the Babylonians." 
Then, on p. 93: ``The Pythagorean school,
which had rigorously established the concept of commensurable quantities,
and attached an almost religious character to it, could not hold to this
point of view; and it is possible that it was the failure of repeated attempts
to express rationally what led them finally to
demonstrate that this number is irrational."

Let me also recall that Bourbaki traces back the concept of group to the Greeks, and, more precisely, to the Pythagoreans. He writes in the chapter \emph{Groups generated by reflections: root systems} of his \emph{Elements of History}:  
``Historically, the beginnings of the theory predate
the introduction of the concept of groups: it actually has its origins
in studies of the ‘regularity’ or ‘symmetries’ of geometric figures,
 and in particular in the determination of regular polygons
and polyhedra (dating back no doubt to the Pythagoreans),
which constitutes the crowning achievement of Euclid's \emph{Elements}, and one of the
most admirable creations of Greek genius."

 In closing this section on Bourbaki, I would like to return to Simone Weil, who, as we have seen, participated to several Bourbaki congresses. In her letter to her brother written in February 1940 which we already quoted  \cite[p. 437]{SWeil-Corresp}, she talks not only about abstract mathematics, but also about concrete problems related to physics. She writes: 
 \begin{quote}\small
 ``I remember that in Chan\c cay or Dieulefit, you said that these studies on Egypt and Babylon raise doubts as to the role of the Greeks as creators in the field of mathematics, a role that had been attributed to them until now. But I believe that, so far (subject to further discoveries), we find confirmation of this role. The Babylonians seem to have focused on abstract exercises involving numbers, while the Egyptians appear to have proceeded in a purely empirical manner. The application of a rational method to concrete problems and to the study of nature seems to have been unique to the Greeks. The odd thing is that the Greeks must have known Babylonian algebra, and yet we find no trace of it among them before Diophantus; for the algebraic geometry of the Pythagoreans is something quite different. There must be religious concepts underlying this; apparently, the secret religion of the Pythagoreans had to accommodate geometry and not algebra.   If the Roman Empire had not destroyed all esoteric cults, we might perhaps understand something about these enigmas.
 \end{quote}

Laurent Lafforgue, in his paper \cite{Lafforgue}, rightly emphasizes that neither Simone Weil nor her brother had much regard for algebra used as an automatism. He recalls that in one of her notebooks, Simone Weil wrote, in capital letters, a succinct formula: ``MONEY, MACHINERY, ALGEBRA. The three monsters of modern civilization. Complete analogy. [\ldots] The institution of algebra corresponds to a fundamental error concerning the human mind" \cite[t. VI, vol; 1, p. 100]{Weil-S}.
We quoted in \S \ref{s:Correspondence} Simone Weil's letter to her brother in which she writes:  ``I believe that God, according to Pythagorean teaching, is a perpetual geometer---but not an algebraist." Echoes of all this can be found in Bourbaki's 
\emph{Elements of the History of Mathematics} \cite[p. 35]{Bourbaki-H}:  

\begin{quote}\small
 If we are to accept as authentic the traditional
maxim ``All things are numbers" of the early Pythagoreans,
we can consider it as evidence of an early attempt
to reduce the geometry and algebra of the time to arithmetic.
Although the discovery of irrational numbers seemed to close
this path forever, the reaction it triggered in Greek mathematics
was a second attempt for a synthesis, this time taking
geometry as a basis and incorporating, among other things, the methods of
solving algebraic equations inherited from the Babylonians [\ldots]
\end{quote}

Such texts could also well be echoes of André Weil's ideas, reformulated by his sister and returned to him.

\bigskip
\bigskip

 \noindent {\bf Acknowledgement.} I would like to warmly thank Bharath Sriraman who read a first version of this paper and kindly sent me several useful comments which helped improving it.

\bigskip

\noindent Athanase Papadopoulos,  Université de  Strasbourg et  Centre National de la Recherche Scientifique.
Institut de Recherche Mathématique Avancée and Centre de Recherche et d'Expérimentation sur l'Acte Artistique, 
7 rue René Descartes,
67084 Strasbourg Cedex France.
Email:  papadop@math.unistra.fr

\end{document}